\documentclass[]{amsart}

\usepackage{lineno,hyperref}
\usepackage{amsmath}
\usepackage{amsfonts}
\usepackage{latexsym}
\usepackage{amsthm}
\usepackage{amssymb, ifthen, enumerate}
\usepackage{color}
\usepackage{ulem}
\usepackage{longtable}
\usepackage{german}
\usepackage{graphicx}
\modulolinenumbers[5]

\begin{document}

\noindent
\Large{\textbf{Vom Gruppentheoriekolloquium 1954 bis zu \glqq B\"aren mit Kindern und Kegeln{\grqq}}}

\noindent
\normalsize Eine Einordnung der fr\"uhen Gruppentheorietagungen am MFO

\noindent
Rebecca Waldecker


\vspace{2cm}

\begin{figure}[h]
\centering
\includegraphics{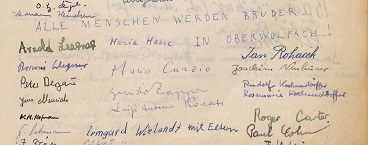}\caption{Auszug aus dem G\"astebuch 1959}

\end{figure}


\noindent
Das Jahr 1954 markiert den Beginn einer Tagungstradition im Bereich der endlichen Gruppentheorie, die bis heute andauert und deren Anf\"ange wir hier beleuchten werden. Dabei geben folgende Fragen Orientierung:

\begin{itemize}
\item
Was f\"ur Tagungen wurden geplant und von wem?

\item
Wie kam es zur thematischen Ausrichtung?

\item
Wer wurde eingeladen, wer hat vorgetragen und wor\"uber?

\item
Was waren die wichtigsten Themen und wer waren die zentralen Personen?

\item
Welche Netzwerke sind erkennbar und wie hat sich die beginnende Tagungstradition weiterentwickelt und etabliert?

\item
Welche Frauen waren zu der Zeit im Bereich der Gruppentheorie aktiv und inwiefern waren sie sichtbar auf diesen Tagungen?

\end{itemize}
\medskip
Als Quellen dienen dabei die G\"asteb\"ucher, Vortragsb\"ucher und Tagungsberichte, soweit vorhanden (siehe\cite{ODA}), Publikationen der Beteiligten aus der Zeit (\cite{MSN}), Korrespondenz, Informationen aus Universit\"atsarchiven und von verschiedenen 
Internetseiten (z.B. \cite{MacTutor},\cite{IMN},\cite{Abel} and \cite{IMU}), sowie Berichte von Zeitzeugen. 
Bei den Photos habe ich, wenn nicht anders vermerkt, auf das Archiv des MFO zur\"uckgegriffen, und die anderen Abbildungen sind selbst erstellte Ausschnitte aus einem G\"astebuch bzw. aus einem Artikel.
Der hier betrachtete Zeitraum wurde eingegrenzt auf ca. 1954 bis 1962, weil in den Jahren 1961/1962 deutlich zu erkennen ist, dass sich die Anzahl und Diversit\"at der Tagungen erh\"oht und dass neue, stark personenbezogene Traditionen beginnen. Diese werden wir zwar am Ende kurz beleuchten, als Ausblick, aber sie verdienen eine separate n\"ahere Betrachtung. Zum Schluss m\"ochte ich noch meine Perspektive einordnen, quasi die \glqq pers\"onliche Brille{\grqq}, durch die ich auf das umfangreiche Material schaue. Zum einen bin ich Mathematikerin und keine Mathematikhistorikerin, und daher gelingt es mir nicht immer, Quellen kritisch einzuordnen. Zum anderen führt meine Forschungstätigkeit im Bereich der endlichen Gruppen  
dazu, dass ich nicht alle Tagungsinhalte in gleicher Weise fachkundig einordnen kann. Mein Schwerpunkt liegt daher auf Themen, die mit Entwicklungen hin zur Klassifikation der endlichen einfachen Gruppen zu tun haben oder z. B. mit Eigenschaften aufl\"osbarer Gruppen oder Permutationsgruppen.

\section{Um welche Tagungen geht es?}

Wir beginnen mit einem Auszug aus einem Bericht von Wilhelm S\"uss (\cite{Suss}). Demnach haben in den Jahren 1951 und 1952 jeweils einmal Tagungen zu den Themen \glqq Moderne Algebra und Zahlentheorie{\grqq} und \glqq Geometrie{\grqq} stattgefunden. Diese Passage bezieht sich auf das erste Thema:
\textit{\glqq Die moderne Algebra und Zahlentheorie, deren urspr\"ungliche Heimat im wesentlichen Deutschland gewesen ist, w\"ahrend sie heute mit außerordentlichem Erfolg an vielen Orten besonders in USA gepflegt wird, hat zweimal in der Berichtszeit eine gr\"oßere Anzahl von Interessenten in Oberwolfach zusammengef\"uhrt. Von den teilnehmenden Ausl\"andern seien genannt R.BAER (Illinois), B.H.NEUMANN (Manchester) und Hanna NEUMANN (Hull). Die bedeutendsten deutschen Vertreter der Algebra waren HASSE und WITT (Hamburg) und Sch\"uler, darunter M.KNESER und P.ROQUETTE. Themen: Erweiterung der Klassenk\"orpertheorie, Funktionenk\"orper, moderne Gruppentheorie.{\grqq}}
Es ist irritierend, wie selbstverst\"andlich Deutschland als Heimat der Gebiete Algebra und Zahlentheorie bezeichnet wird und dass ausgerechnet
drei Deutsche, die w\"ahrend der Nazi-Diktatur ins Ausland geflohen waren, hier als \glqq Ausl\"ander{\grqq} aufgef\"uhrt werden und damit dem MFO eine gewisse Internationalit\"at bescheinigen sollen\footnote{In \cite{Remmert} wird dies ab Seite 260 detailliert beleuchtet.}. 
Beim Thema Algebra w\"aren zu der Zeit noch die Namen Wolfgang Gasch\"utz und Helmut Wielandt erwartbar gewesen, ebenso Ruth Moufang beim Thema Geometrie, aber in der Namensliste am Ende des Berichts tauchen sie nicht auf. Bartel L. van der Waerden und Olga Taussky-Todd dagegen stehen in der Namensliste, werden aber von S\"uss nicht explizit genannt. Ab 1954 finden 
regelm\"a{\ss}ig Tagungen zu Themen aus der Gruppentheorie statt, was ich auf die im Bericht erw\"ahnte \textit{\glqq gr\"oßere Anzahl an Interessenten{\grqq}} zur\"uckf\"uhre und das gut dokumentierte Interesse von S\"uss, sichtbare Aktivit\"aten am MFO zu erzeugen und damit die Bedeutung und Au{\ss}enwirkung des Instituts zu steigern. Gemeinsam gelingt es Reinhold Baer und Wilhelm S\"uss, innerhalb weniger Jahre eine Tagungstradition zu etablieren mit zun\"achst einer Konferenz pro Jahr (mit Ausnahme von 1957) und ab 1961 sogar mehreren Konferenzen j\"ahrlich mit unterschiedlichen Zielgruppen. Im Zentrum dieses Artikels stehen die dokumentierten Tagungen aus den Jahren 1954 bis 1962:

\begin{itemize}
\item
Kolloquium \"uber Gruppentheorie und Grundlagen der Geometrie, 7. -- 11. Juni 1954, geleitet von Reinhold Baer (Workshop 5423).
\item
Gruppentheorie 2. -- 15. August 1955, geleitet von Helmut Wielandt (Workshop 5531).
\item
Gruppentheorie 25. -- 29. September 1956, geleitet von Helmut Wielandt (Workshop 5639).
\item
Die Geometrien und ihre Gruppen 26. -- 31. Mai 1958, geleitet von Reinhold Baer und Jacques Tits (Workshop 5822).
\item
Tagung Gruppentheorie 19. -- 22. Mai 1959, geleitet von Reinhold Baer und Helmut Wielandt (Workshop 5921).
\item
Kolloquium \"uber Gruppentheorie 6. -- 12. Juni 1960, geleitet von Reinhold Baer und Helmut Wielandt (Workshop 6023)
\item
Kolloquium \"uber abelsche Gruppen, 5. -- 10. M\"arz 1961, geleitet vermutlich von Reinhold Baer (Workshop 6110).
\item
Gruppentheorietagung 16. -- 21. Oktober 1961, geleitet von Reinhold Baer und Helmut Wielandt (Workshop 6142).
\item
B\"aren mit Kindern und Kegeln, Januar 1962.
\item
Endliche Geometrien und ihre Gruppen, 1. -- 8. April 1962, geleitet von Peter Dembowski (Workshop 6214). Gro{\ss}e personelle Überschneidung mit der Geometrie-Tagung im Juni 1962.
\item
Arbeitstagung Baer, Juni 1962.
\item
Gruppentheorie, 5. -- 11. August 1962, geleitet von Reinhold Baer und Helmut Wielandt (Workshop 6232).
\end{itemize}

Die L\"ucke im Jahr 1957 kann damit erkl\"art werden, dass es in dem Jahr in T\"ubingen ein Kolloquium der Internationalen Mathematischen Union \"uber endliche Gruppen gab (\cite{IMU}). Typisch ist eine Mischung aus gruppentheoretischen und geometrischen Themen und insbesondere endlicher und unendlicher Gruppentheorie, was sich bis heute im Tagungsformat \glqq Gruppen und Geometrien{\grqq} erhalten hat. Zur Bandbreite geh\"oren die gerade im Entstehen begriffene Computeralgebra, die bereits etablierte, aber noch nicht so stark ausdifferenzierte Darstellungstheorie, und auch die algebraischen Gruppen. Erst deutlich sp\"ater gibt es in diesen Bereichen eigene, spezialisierte Tagungsformate. 

Bevor wir uns n\"aher mit der Tagungsorganisation befassen, werfen wir einen Blick auf den durch die Berichte und Briefwechsel auff\"alligen Personenkreis, denn hier zeigen sich Schl\"usselfiguren f\"ur diese fr\"uhe Phase der Oberwolfachtagungen. Das sind die Personen, die in der Forschung sehr aktiv und sichtbar, teils auch untereinander stark vernetzt waren, die Tagungen organisiert und geleitet haben bzw. die hinter den Kulissen in enger Korrespondenz miteinander und mit der Leitung des MFO standen (Wilhelm S\"uss, sp\"ater Hellmuth Kneser, dann Theodor Schneider). H\"aufig waren das auch \glqq Dauerg\"aste{\grqq} am MFO: Reinhold Baer, Wolfgang Gasch\"utz, Bertram Huppert, Bernhard und Hanna Neumann, Jacques Tits und Helmut Wielandt. Schauen wir sie uns in alphabetischer Reihenfolge etwas genauer an.


\section{Schl\"usselfiguren}

\textbf{Reinhold Baer} (1902-1979), ein Sch\"uler von Hellmuth Kneser, wird der Noether-Schule zugerechnet (siehe \cite{Koreuber}), verbrachte aber vor der Promotion auch ein Jahr in Kiel (1924, mit einem Stipendium) bei Helmut Hasse u.a. und war 1926-28 Assistent von Alfred Loewy in Freiburg.

\begin{figure}[h]
\centering
\includegraphics[scale=0.5]{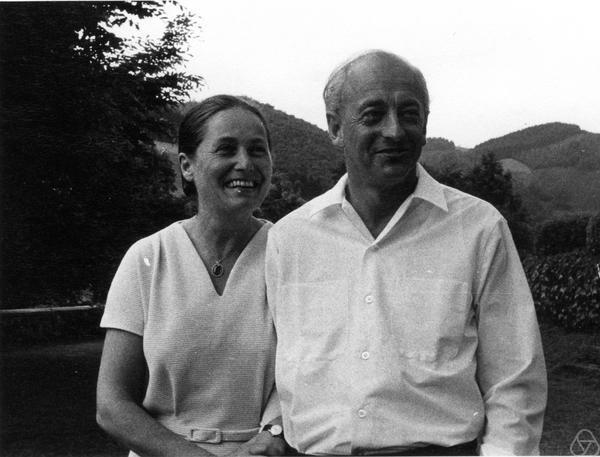}\caption{Marianne und Reinhold Baer}

\end{figure}

Zu der Zeit arbeiteten mehrere Personen in Freiburg, die stark von Emmy Noether beeinflusst waren\footnote{Siehe dazu auch \cite{Remmert1}.} – neben den genannten auch Wolfgang Krull. Nach der Habilitation 1928 in Freiburg wechselte Baer nach Halle an der Saale, wo ihm allerdings 1933 die Lehrerlaubnis entzogen wurde. Nach Stationen in Manchester, Oxford und Princeton blieb er ab 1938 f\"ur l\"angere Zeit an der University of Illinois, Urbana (USA). Von dort aus nahm er an der Tagung im Jahr 1954 teil, die unter seiner Leitung stattfand, und er trug dort \"uber Noethersche \"uberaufl\"osbare Gruppen vor. 1956 hielt Baer einen Vortrag \"uber Engelsche Elemente\footnote{Dazu gibt es eine bis 2020 viel zitierte Publikation, siehe \cite{Baer1957}.}, wobei er nicht nur endliche Gruppen untersucht hatte. Er bewies u.a. folgende Charakterisierung: Eine Noethersche Gruppe ist dann und nur dann von endlicher Klasse (d.h. die absteigende Kommutatorreihe endet mit endlichem Index), wenn sie von ihren Engelschen Elementen erzeugt wird. Weitere Themen von Baer: Sylowturmgruppen (in einem Artikel 1955 konzeptioniert durch Überaufl\"osbarkeit, siehe auch \cite{Baer1958}), Partitionen, und sog. Hauptuntergruppen. Ein wichtiges Resultat, von dem es zahlreiche Variationen gibt, ist folgendes: \textit{\glqq A finite group G is a Sylow tower group if and only if every pair of elements of G generates a Sylow tower group.{\grqq}}

\"Uberaufl\"osbarkeit etablierte sich sp\"ater als Begriff zwischen Aufl\"osbarkeit und Nilpotenz in dem Sinne, dass jede nilpotente Gruppe \"uberaufl\"osbar ist (aber nicht umgekehrt) und dass jede \"uberaufl\"osbare Gruppe aufl\"osbar ist (aber nicht umgekehrt). Wir haben bereits gesehen, dass Reinhold Baer mit einigen Kollegen in Deutschland gut vernetzt war und gemeinsam mit Bernhard Neumann stark in die Planung der Tagungen involviert war, etwa, was die Schwerpunktthemen und die Einzuladenden betrifft. Daher spielte er bei der Tagungsorganisation auch schon vor seiner R\"uckkehr nach Deutschland (Professur in Frankfurt a.M. ab Juli 1956) eine wichtige Rolle. Hintergr\"unde dazu findet man in einigen Arbeiten von Volker Remmert\footnote{etwa \cite{Remmert} und \cite{Remmert2}}. Hinweise auf den gro{\ss}en Einfluss Baers auf die Wissenschaftslandschaft in mehreren Bereichen der Mathematik liefern neben seinem mathematischen Werk\footnote{MathSciNet (\cite{MSN}) listet am 9.6.2021 150 Publikationen und 1044 Zitate.}
seine zahlreichen Sch\"uler*innen\footnote{Maths Genealogy (\cite{MGP}) am 28.12.2020: \glqq Reinhold Baer has 60 students and 937 descendants{\grqq}.}, wobei die meisten mathematischen Nachkommen \"uber Paul Conrad und Gerhard Michler zustande kommen. 

\vspace{1cm}
\textbf{Wolfgang Gasch\"utz} (1920-2016), ein Sch\"uler von Karl-Heinrich Weise und in Kiel t\"atig, arbeitete in dem hier behandelten Zeitraum u.a. an Gruppen, in denen das Normalteilersein transitiv ist, zu aufl\"osbaren Gruppen und zur Eulerschen Funktion (\cite{Gaschutz1959}). 1954 trug er vor zu \textit{\glqq einem einfacheren Beweis einer bereits bekannten Aussage{\grqq}} \"uber treue irreduzible Darstellungen endlicher Gruppen. Er war bei der Tagung 1955 nur kurz dabei -- eine erbetene Verschiebung klappte nicht. \footnote{Briefwechsel S\"uss-Gasch\"utz im Juni und Juli 1955, siehe Universit\"atsarchiv Freiburg (UAF) C89 / 297 ab S. 19.} Konkret schrieb er am 10.7.1955 an Wilhelm S\"uss\footnote{UAF C89 / 297 S. 26}: \textit{\glqq Ich habe mich nun doch entschlossen, daran teilzunehmen und werde etwa am 6. oder 7. August allein dort eintreffen. Falls es das Programm gestattet, kann ich zwei Vortr\"age halten: 1) Erzeugendenzahl und Existenz von p-Faktorgruppen, 2) Gruppen, deren minmale Untergruppen Normalteiler sind, je Thema w\"urde ich etwa 1 Stunde ben\"otigen.{\grqq}}  Es gibt eine Unterstreichung und ein Pluszeichen am Rand der Notiz, und das Vortragsbuch zeigt ebenfalls, dass er die beiden vorgeschlagenen Vortr\"age gehalten hat. 
Das zweite genannte Vortragsthema f\"uhrte zur Publikation \cite{Gaschutz1957}, in der er schreibt: \textit{\glqq Die Vielgestaltigkeit der endlichen Gruppentypen beruht zum Teil darauf, da{\ss} f\"ur sie nicht allgemein die folgende Transitivit\"at gilt: Ein Normalteiler eines Normalteilers ist ein Normalteiler.{\grqq}}
 Falls die Transitivit\"at doch gilt, nennt er die Gruppe t-Gruppe (Definition). Das basiert auf Ideen von Ernest Best und Olga Taussky (\cite{BT1942}), und es gab schon Resultate zum Spezialfall, dass alle Sylow-Untergruppen zyklisch sind. Im Artikel konzentriert sich Gasch\"utz auf aufl\"osbare t-Gruppen und beschreibt ihre Struktur, erarbeitet Charakterisierungen und verallgemeinert dann Resultate auf allgemeine Gruppen \"uber die Betrachtung des aufl\"osbaren Radikals. Dies ist seine am h\"aufigsten zitierte Arbeit, bis ins Jahr 2019 sind Referenzen darauf zu finden (\cite{MSN}). Gasch\"utz war im hier betrachteten Zeitraum sehr pr\"asent und korrespondierte viel mit S\"uss, wie wir noch sehen werden, und das oft in einem herzlichen, vertrauten Ton. 1960 hielt Gasch\"utz einen Vortrag zur Eulerschen Funktion. Dabei ist f\"ur jede Gruppe G und jede nat\"urliche Zahl s definiert: $\varphi_G(s)$ ist die Anzahl der s-elementigen Erzeugendensysteme von G. Die Definition geht auf Philipp Hall zur\"uck. Gasch\"utz bewies Resultate dazu, wie sich die $\varphi$-Werte von Normalteilern und Faktorgruppen verhalten in Bezug auf die $\varphi$-Werte der Ausgangsgruppe, was Reduktionsschritte m\"oglich macht. Das Thema ist u.a. in der Computeralgebra relevant und f\"ur probabilistische Methoden. Im ersten Zitat dazu aus dem Jahr 1998 geht es zum Beispiel um Algorithmen f\"ur endliche Gruppen und speziell um Markov-Ketten auf der Menge der Erzeugendensysteme einer gegebenen Gr\"o{\ss}e, f\"ur eine feste Gruppe (\cite{Diaconis1998}). Gasch\"utz‘ Arbeit wird u.a. deshalb bis in die Gegenwart hinein zitiert\footnote{MathSciNet (\cite{MSN}) listet insgesamt 35 Publikationen und 533 Zitate.}. Einer seiner mathematischen Nachkommen\footnote{Maths Genealogy (\cite{MGP}) am 28.12.2020: \glqq Wolfgang Gasch\"utz has 22 students and 136 descendants{\grqq}.}, Joachim Neub\"user, wird uns sp\"ater noch begegnen.

\textbf{Bertram Huppert} (geb. 1927) ist ein Sch\"uler von Helmut Wielandt. Nach dem Studium in Mainz, bei Wielandt, folgte er ihm 1951 nach T\"ubingen und wurde dort 1953 als erster Doktorand Wielandts promoviert. Er geh\"ort zu den Vielvortragenden, au{\ss}erdem sind die Vortragsthemen meist sehr klar seinen Publikationen zuzuordnen\footnote{Siehe dazu \cite{Hupp1953} -- \cite{Hupp1962}.}. 1954 ging es um Aufl\"osbarkeitss\"atze f\"ur faktorisierbare Gruppen, 1955 um Hauptreihen aufl\"osbarer Gruppen und zweifach transitive aufl\"osbare Permutationsgruppen. 1956 trug er vor \"uber Normalteiler mehrfach transitiver Gruppen und \"uber endliche Engel-Gruppen, dann wieder 1960 \"uber subnormale Untergruppen und p-Sylowgruppen. Das Konzept von subnormalen Untergruppen, heute Subnormalteiler genannt, war damals neu und wurde von Helmut Wielandt als \glqq nachinvariante Untergruppen{\grqq} eingef\"uhrt (Siehe dazu \cite{Wielandt1939}.). Schlie{\ss}lich trug Bertram Huppert 1961 \"uber scharf dreifach transitive Permutationsgruppen vor. Er wird im weiteren Verlauf der Tagungsreihen noch eine Rolle spielen und ist au{\ss}erdem durch seine zahlreichen B\"ucher bekannt\footnote{Vor allem die Lehrb\"ucher \cite{HuppEG1} -- \cite{HuppChar}.}. MathSciNet listet 57 Publikationen und 4399 Zitate, von denen fast 4000 auf das Konto seiner oben genannten B\"ucher gehen. Allein f\"ur sein Buch \glqq Endliche Gruppen I{\grqq} sind \"uber 2800 Zitationen vermerkt. Bemerkenswert ist Hupperts Rolle, als sich die Tagungen in Oberwolfach inhaltlich st\"arker ausdifferenzieren – sowohl f\"ur das Thema Gruppentheorie als auch f\"ur das Thema Darstellungstheorie ist er sp\"ater als Tagungsorganisator t\"atig. 

\begin{figure}[h]
\centering
\includegraphics[scale=0.8]{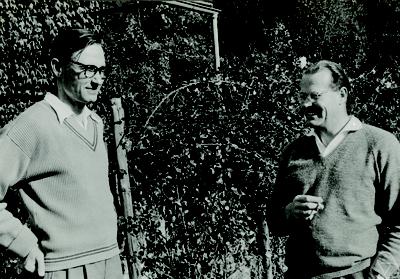}\caption{Huppert und Gasch\"utz 1961}

\end{figure}

\textbf{Bernhard Neumann} (1909-2002) war ebenfalls mehrmals bei den Tagungen dabei und trug dort vor. Er war ein Sch\"uler von Issai Schur in Berlin, schrieb aber nach seiner Flucht aus Deutschland (im April 1933) eine weitere Arbeit bei Philipp Hall in Cambridge, UK, und war dann l\"angere Zeit in Manchester t\"atig. Wie sich im n\"achsten Abschnitt zeigen wird, war er oft an der Tagungsorganisation beteiligt, an Terminabsprachen, er machte Vorschl\"age f\"ur Einzuladende und wirkte als Verbindungsperson zu weiteren Kolleg*innen in Gro{\ss}britannien. Zusammen mit Reinhold Baer und Friedrich Wilhelm Levi, dessen Namen wir bei diesen fr\"uhen Tagungen auch mehrmals in den G\"asteb\"uchern finden, geh\"ort er zu einem Kreis von Personen, um deren Anwesenheit am MFO sich Wilhelm S\"uss nachweislich sehr bem\"uhte. Volker Remmert schreibt dazu\footnote{Seite 362 in \cite{Remmert}.}: \textit{\glqq Neumann organized his first conference on group theory at Oberwolfach in 1955. ... Both Baer and Neumann were instrumental to the remigration of
mathematical ideas and theories to Germany, in particular in group theory.{\grqq}} Wir kommen darauf im Abschnitt 
\"uber Netzwerke noch einmal zur\"uck -- S\"uss verfolgte ganz klar das Ziel, das MFO und seine Aktivit\"aten sichtbarer und einflussreicher zu machen und hatte erkannt, dass gute Kontakte zu Schl\"usselfiguren in den entsprechenden Fachcommunities dabei wichtig sind.
 1955 hielt Bernhard Neumann zwei Vortr\"age: \glqq Aufsteigende Reihen von Kommutatorgruppen{\grqq} und \glqq Ordnung von Automorphismengruppen endlicher Gruppen{\grqq}. Zum zweiten Thema erschienen sp\"ater auch Publikationen, die oft zitiert wurden (siehe zum Beispiel \cite{Ledermann}). Die Frage nach der Anzahl von Automorphismen einer Gruppe geht zur\"uck auf Garrett Birkhoff und Philip Hall (siehe \cite{BH}), und Neumann beweist eine Vermutung zu Untergrenzen der Anzahl ab einer gewissen Gruppengr\"o{\ss}e. Er vermutet, dass die Grenze nicht bestm\"oglich ist und motiviert das mit abelschen Gruppen. Au{\ss}erdem ist bei unendlichen Gruppen alles anders! Etwa kann eine unendliche Gruppe auch nur wenige, z.B. nur zwei, Automorphismen haben. Interessante Spezialf\"alle sind die eines fixpunktfreien Automorphismus der Ordnung 2 (seit 1940 bekannt, siehe \cite{NeumannB1956}) oder 3, besonders der Zusammenhang zur Nilpotenz. Das Thema blieb aktuell, etwas sp\"ater bewies John G. Thompson in seiner Doktorarbeit ein bis heute viel verwendetes Resultat dazu, das jetzt zum Standardrepertoire der endlichen Gruppentheorie geh\"ort:
\textit{Eine endliche Gruppe, die einen fixpunktfreien Automorphismus von Primzahlordnung hat, ist nilpotent.} Siehe \cite{Thompson}. 
Hinweise auf Bernhard Neumanns Einfluss sehen wir nicht nur in seinem umfangreichen mathematischen Werk \footnote{MathSciNet (\cite{MSN}) listet f\"ur Neumann 125 Publikationen und 1270 Zitate.}, sondern auch bei der gro{\ss}en Zahl an Nachkommen in verschiedenen Bereichen der Mathematik \footnote{Maths Genealogy am 28.12.2020: \glqq Bernhard Neumann has 18 students and 351 descendants{\grqq}.}.

\begin{figure}[h]
\centering
\includegraphics[scale=0.8]{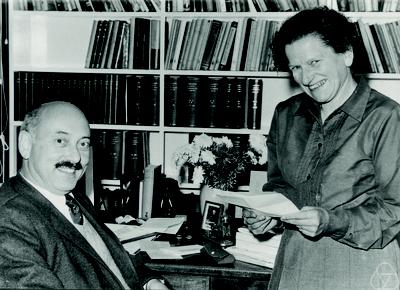}\caption{Bernhard und Hanna Neumann}

\end{figure}

\textbf{Hanna Neumann} (1914-1971), geboren als Johanna von Caemmerer, war eine Sch\"ulerin von Olga Taussky-Todd und zur Zeit ihrer Besuche in Oberwolfach in Hull, UK, sp\"ater auch in Manchester t\"atig. Studiert hatte sie an der Universit\"at Berlin, wo sie u.a. Helmut Wielandt und Bernhard Neumann kennenlernte, und ihr Werdegang ist stark von der politischen Situation gekennzeichnet\footnote{Siehe zum Beispiel \cite{MacTutor}.}. Das 1932 begonnene Studium schloss sie in Berlin 1936 mit einem Staatsexamen ab und wechselte zur Promotion nach G\"ottingen, zu Helmut Hasse. Bernhard Neumann, mit dem sie seit 1933 mindestens freundschaftlich verbunden war, hatte Deutschland bereits 1933 verlassen, und 1938 beschloss Hanna Neumann, die geplante Promotion in G\"ottingen abzubrechen und Bernhard Neumann nach Cardiff zu folgen. Dort heirateten sie heimlich, im Dezember 1938. Hanna Neumann war eine der wenigen Frauen, die in dieser Zeit mehrmals bei den Tagungen in Oberwolfach dabei waren und selbst zahlreiche mathematische Nachkommen\footnote{Maths Genealogy (\cite{MGP}) am 25.5.2021: \glqq Hanna Neumann has 10 students and 68 descendants.{\grqq}.} hatten. Sie trug 1956 vor zum Durchschnitt von Untergruppen freier Gruppen, zu einer Erg\"anzung zu Jordan-H\"older und zu Normalreihen (leider im Vortragsbuch nicht vollst\"andig lesbar). 1960 war sie dabei mit dem Vortrag \glqq Linked products of groups{\grqq}, siehe dazu \cite{NeumannH1960a} und \cite{NeumannH1960b}. MathSciNet (\cite{MSN}) listet f\"ur Hanna Neumann 30 Ver\"offentlichungen und 673 Zitate, wovon \"uber 300 zum Buch \"uber Gruppenvariet\"aten geh\"oren und \"uber 150 zur Arbeit \"uber Einbettungss\"atze (siehe \cite{NeumannH1960a}).

\begin{figure}[h]
\centering
\includegraphics{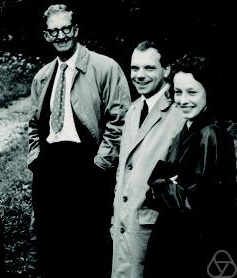}\caption{M. Lazard, Jacques Tits und Marie-Jeanne Tits 1957}

\end{figure}

\textbf{Jacques Tits} (1930-2021), ein Sch\"uler von Paul Libois, spielte eine gro{\ss}e Rolle u.a. wegen seiner geometrischen Sichtweise auf bereits bekannte Gruppen, etwa Suzukigruppen (siehe \cite{Tits}), aber auch bei der Konstruktion neuer Gruppen. Er wirkte zu der Zeit in Br\"ussel und war nicht nur bei mehreren Gruppentheorietagungen dabei, sondern auch regelm\"a{\ss}ig Gast auf Geometrie-Tagungen, und im hier behandelten Zeitraum schrieb er grunds\"atzlich in franz\"osischer Sprache ins Vortragsbuch. 
Tits trug 1955 vor \"uber mehrfach transitive Transformationsgruppen und \"uber Geometrien, die zu halb-einfachen Lie-Gruppen geh\"oren, dann 1958 \"uber eine geometrische Deutung der Coxeter-Weyl-Witt-Dynkinschen Schemata und 1960 \"uber die Definition gewisser einfacher algebraischer Gruppen verm\"oge Erzeugern und Relationen und \"uberaufl\"osbare Lie-Gruppen. 
MathSciNet listet 161 Ver\"offentlichungen und die beeindruckende Zahl von 5783 Zitaten, von denen viele auf seine Arbeiten \"uber reduktive Gruppen zur\"uckgehen. Er gilt als Begr\"under der Theorie der Geb\"aude, damit auch als wegweisend in dem Bereich der Mathematik, der oft als geometrische Gruppentheorie bezeichnet wird. Er wurde vielfach f\"ur seine akademischen Leistungen ausgezeichnet, u.a. 2008 mit dem Abelpreis.\footnote{Auf der Seite \cite{Abel} sind alle bisherigen Preistr\"ager*innen verzeichnet.}

\begin{figure}[h]
\centering
\includegraphics{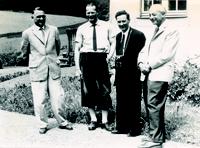}\caption{Helmut Wielandt (links), mit Wolfgang Gasch\"utz, Jen\"o Sz\`ep und Reinhold Baer 1959}

\end{figure}

\textbf{Helmut Wielandt} (1910-2001), ebenfalls ein Sch\"uler von Issai Schur, war in dem hier betrachteten Zeitraum in T\"ubingen t\"atig und hatte h\"aufig internationale G\"aste (z.B. John G. Thompson im Sommer 1958, auch in den Briefwechseln mit Wilhelm S\"uss ist das immer wieder Thema, siehe z.B. \cite{MacTutor}). 
Er sprach bei der Tagung 1954 \"uber faktorisierte Gruppen, 1955 u.a. \"uber Permutationsgruppen von Primzahlgrad, und sp\"ater kamen noch andere Themen dazu, z.B. Subnormalteiler. Weiterhin entwickelte Wielandt Gedanken f\"ur eine allgemeine Theorie der Permutationsgruppen - aus der Zeit bis 1960 stammen viele gute und einflussreiche Arbeiten. Es gab intensive Korrespondenz mit Wilhelm S\"uss, Hellmuth Kneser u.a., \"uber Kontakte zu Reinhold Baer, Tagungsinhalte (1954) und Einzuladende, und er leitete auf Wunsch von S\"uss die Tagung 1955, deren Erfolg wir im n\"achsten Abschnitt beleuchten werden. Dort werden wir auch sehen, dass er sich zwischenzeitlich etwas aus der Organisation zur\"uckgezogen hat, und warum. Er war aber z.B. wieder beteiligt an der Planung des Pfingstkolloquiums 1959, und hier finden wir detaillierte Absprachen im Briefwechsel mit der Institutssekret\"arin Anneliese Pfefferle im Fr\"uhjahr 1959\footnote{UAF C89 / 397 ab S. 93.}. Am 14.4.1959 schrieb er zum Beispiel ausf\"uhrlich drei Listen auf und bat um Verschickung der \textit{\glqq endg\"ultigen Einladung zum Gruppenkolloquium{\grqq}}\footnote{UAF C89 / 397, S. 96 oben.}. Dabei stehen auf Liste 1 diejenigen, die auf die erste Einladung positiv reagiert haben, auf Liste 2 diejenigen, die von Teilnehmern auf Liste 1 vorgeschlagen worden sind, und auf Liste 3 stehen \textit{\glqq die folgenden, die die erste Einladung vielleicht verbummelt haben{\grqq}}\footnote{UAF C89 / 397, S. 96 etwa mittig.}.
Auch zu Absagen gibt es Notizen in dem Brief. In einem Nachtrag am 15.4.1959 wurden noch drei Sch\"uler von Wolfgang Gasch\"utz nachgemeldet, die eine Einladung erhalten sollen, darunter der bereits erw\"ahnte Joachim Neub\"user. Der Andrang auf die Oberwolfach-Tagungen wurde merklich gr\"o{\ss}er, das Publikum internationaler, und gemeinsam mit Reinhold Baer organisierte Wielandt noch mehrere Gruppentheorietagungen in den fr\"uhen 1960er Jahren. Sp\"ater finden sich im ODA (\cite{ODA}) auch Berichte \"uber Tagungen, die Wielandt gemeinsam mit Hans Zassenhaus organisiert hat. 
Wielandts Wirken ist weit \"uber den Bereich der Permutationsgruppentheorie sichtbar\footnote{MathSciNet listet 65 Publikationen und 1261 Zitate, von denen etwa die H\"alfte sich auf sein Standardwerk zu endlichen Permutationsgruppen (\cite{Wielandt}) bezieht.}, und bekannte Sch\"uler\footnote{Maths Genealogy zeigt am 10.3.2020 zu Helmut Wielandt \glqq 22 students and 544 descendants{\grqq}, dabei stammen die meisten mathematischen Nachkommen von Harro Heuser.} sind neben Bertram Huppert auch Harro Heuser und Olaf Tamaschke.\\

Wenn wir uns auf diese Schl\"usselfiguren konzentrieren, so sehen wir in Deutschland die Standorte Freiburg (Wilhelm S\"uss), T\"ubingen (Helmut Wielandt, Bertram Huppert, sp\"ater auch der Pickert-Sch\"uler Helmut R. Salzmann), Kiel (Wolfgang Gasch\"utz) und Frankfurt (Reinhold Baer ab 1956).
Wir kommen auf diese Standorte und die Verbindungen im Abschnitt \"uber Netzwerke zur\"uck.

\section{Was wurde geplant, wie und von wem?}

Ausschnitte aus Briefen und aus einigen Tagungsberichten geben uns einen ersten Eindruck von den Tagungen 1955 und 1956, von der damit zusammenh\"angenden Organisation inkl. Finanzierung\footnote{PS im Brief von Wilhelm S\"uss an Helmut Wielandt am 9.5.1956, UAF C89 / 397 S. 89, auch Thema im Brief von S\"uss an Wolfgang Gasch\"utz am 21.6.1956, UAF C89 / 297 S. 23.}, der Atmosph\"are dort und der thematischen Vielfalt.
S\"uss schreibt an Wielandt am 24.8.55: \textit{\glqq … m\"ochte ich Ihnen nochmals f\"ur den Erfolg der Tagung herzlich danken, der ja weitgehend auch Ihr Werk gewesen ist. Was Sie am letzten Tag vorgetragen haben, entnahm ich aus den begeisterten Worten der Anderen und aus Ihrem Protokoll. Ich hatte den Eindruck, da{\ss} allein schon dieser Fund Ihnen reichlich genug Belohnung f\"ur die gehabte Arbeit sein konnte. Hoffentlich ist es auch so!{\grqq}}\footnote{UAF C89 / 397 S. 72} 

Laut Vortragsbuch gab es insgesamt drei Beitr\"age von Wielandt bei der Tagung 1955, darunter zwei an den letzten beiden Tagen. Im letzten Vortrag ging es unter dem Titel \glqq Ein Beweis des Ergodensatzes“ um eine Gruppenwirkung auf einer abgeschlossenen, beschr\"ankten Teilmenge eines Hilbertraums und um die Existenz eines Fixpunkts. S\"uss k\"onnte das gemeint haben, aber ich konnte keine dazu passende Publikation finden oder andere Hinweise darauf, dass ausgerechnet dieser Beitrag von Wielandt auf besonders viel Resonanz gesto{\ss}en ist. Anders verh\"alt es sich mit dem Thema seines Vortrags am vorletzten Tag \"uber \glqq Permutationsgruppen von Primzahlgrad“ – daf\"ur ist Wielandt bekannt und es handelt sich um eine sch\"one Anwendung von Methoden der Darstellungstheorie auf ein Problem aus der Theorie der Permutationsgruppen\footnote{Siehe \cite{Wielandt1956}.}.
In einem Brief vom 16.3.1956 , der sich bereits mit der Planung der n\"achsten Gruppentheorietagung in Oberwolfach befasst, wird Wielandt erneut von S\"uss gelobt f\"ur die \textit{\glqq so gl\"anzende Gruppentagung“}. Bei den Absprachen spielt neben der Terminfindung (u.a. Abstimmung mit der DMV-Tagung 1956) die Einbeziehung von Bernhard Neumann eine Rolle. Dieser hatte nicht nur Interesse an einer Gruppentheorietagung, sondern ist auch f\"ur die rechtzeitige Ansprache der Kolleg*innen in Gro{\ss}britannien zust\"andig\footnote{Ebenfalls aus dem Brief von S\"uss an Wielandt am 16.3.1956, UAF C89 / 397 S. 85}. Wielandt wiederum \"au{\ss}ert in seinem Brief vom 24.3.1956 an S\"uss\footnote{UAF C89 / 397 S. 86} die Hoffnung, zur Tagung auch Kollegen aus Russland und Ungarn gewinnen zu k\"onnen. Leider liegt zur Tagung 1956 kein Bericht vor, aber dem G\"astebuch und dem Vortragsbuch (siehe \cite{ODA}) l\"asst sich entnehmen, dass es gelungen ist, im September einen passenden Termin f\"ur die Tagung zu finden und mehrere Kolleg*innen aus dem Ausland dabei zu haben (laut Vortragsbuch Harald Bergstr\"om, Frans Loonstra, Graham Higman, G. L. Britton, Hanna und Bernhard Neumann, J. A. H. Shepperd, D. H. McLain, John (Sean) Tobin, L\'aszl\'o Fuchs). Die Einladung russischer Kolleg*innen scheiterte allerdings. Als Begr\"undung wird die Zeitplanung angegeben, unter Berufung auf eine vertrauliche Mitteilung des Kultusministeriums Baden-W\"uettemberg\footnote{UAF C89 / 297 S. 25}, denn es wird als unm\"oglich eingesch\"atzt, rechtzeitig \textit{\glqq mit dem Ausw\"artigen Amt F\"uhlung zu nehmen{\grqq}}.
Aus anderen Briefen\footnote{S\"uss an Wielandt UAF C89 / 397 S. 85 und S\"uss an Alexander Ostrowski UAF C89 / 343 S. 66} geht hervor, dass Bernhard Neumann in beiden Jahren, 1955 und 1956, die Veranstaltung einer Tagung angeregt hatte, dass diese aber schlie{\ss}lich haupts\"achlich von Wielandt organisiert wurden. Zumindest im Jahr 1956 war das durchaus problematisch f\"ur Wielandt. Am 25.4.1956 entschuldigt er sich f\"ur seine sp\"ate Antwort\footnote{UAF C89 / 397 S. 88} und schreibt u.a.:
\textit{\glqq Es geht mir im Augenblick ziemlich schlecht, weil das frisch \"ubernommene Dekanat mich wegen Erkrankung der Fakult\"atssekret\"arin \"uberm\"a{\ss}ig viel Kraft kostet. Ich habe, um Verz\"ogerungen zu vermeiden, Herrn Neumann gebeten, die Vorbereitungen mit Ihnen zusammen vorerst allein in Gang zu bringen, bis ich wieder Luft habe.{\grqq}}
Bezugnehmend auf den R\"uckzug von Wielandt schreibt Gasch\"utz am 12.6.1956 an S\"uss:
\textit{\glqq Ich m\"ochte aufgrund dieser Mitteilungen empfehlen, dass nach dem Ausfall von Prof. Wielandt m\"oglichst Herr Prof. Baer die Protektion f\"ur die Tagung \"ubernimmt, falls er dazu bereit ist; er soll im Juli in Deutschland eintreffen.{\grqq}}\footnote{UAF C89 / 297 S. 25}
Dabei bezieht sich Gasch\"utz darauf, dass Baer 1933 aus Nazi-Deutschland geflohen war. Nach mehreren akademischen Stationen war Baer seit 1938 an der University of Illinois (Urbana, USA) t\"atig und hatte dort ma{\ss}geblichen Anteil daran, Illinois zu einem Zentrum der Forschung zu endlichen Gruppen zu machen. Gleichzeitig ist sein Interesse, nach Deutschland zur\"uckzukehren, gut dokumentiert, u.a. in Briefwechseln mit S\"uss (siehe auch \cite{MacTutor} und \cite{Remmert2}\footnote{Ab Seite 36 dort.}). 
Im gleichen Brief von Gasch\"utz an S\"uss wird die Entscheidung f\"ur eine eher \glqq unendliche{\grqq} Gruppentheorietagung mit dem 1957 in Deutschland stattfindenden IMU-Kolloquium zu endlichen Gruppen begr\"undet (siehe \cite{IMU}). 
Zum Schluss bekr\"aftigt Gasch\"utz, dass Baer bei der Organisation eine tragende Rolle spielen soll:
\textit{\glqq Auf jeden Fall erscheint es mir dem Vorhaben dienlicher, wenn ein repr\"asentativer Herr die Einladungen vollzieht. Ich bin gern bereit im Einvernehmen mit dem Betreffenden die technische Seite der Sache in die Hand zu nehmen.{\grqq}}

Genau so kommt es – wie wir schon der Tagungsliste oben entnehmen k\"onnen, ist Reinhold Baer ab 1958 fester Bestandteil des Organisationsteams f\"ur die Oberwolfachtagungen im Bereich der Gruppentheorie. 
Im Jahr 1958 organisiert er die Tagung \glqq Die Geometrien und ihre Gruppen{\grqq} gemeinsam mit Jacques Tits. In Briefwechseln mit S\"uss liest es sich so, als sei Tits bereits f\"ur die Gruppentheorie-Tagung 1955 zur Mitorganisation und nicht nur zur Teilnahme angefragt worden. Im Jahr 1958 l\"auft dann viel Korrespondenz \"uber S\"uss, der zwischen Baer und Tits Vorschl\"age f\"ur Teilnehmer*innen hin- und herschickt und den Stand der Planung im Auge beh\"alt\footnote{UAF C 89/385, S. 49-55}
Ein Auszug aus dem Bericht gibt uns einen Eindruck von der Tagung\footnote{Siehe \cite{ODA}, Bericht zum Workshop 5822.}:

\begin{figure}[ht]
\centering
\includegraphics[scale=0.65]{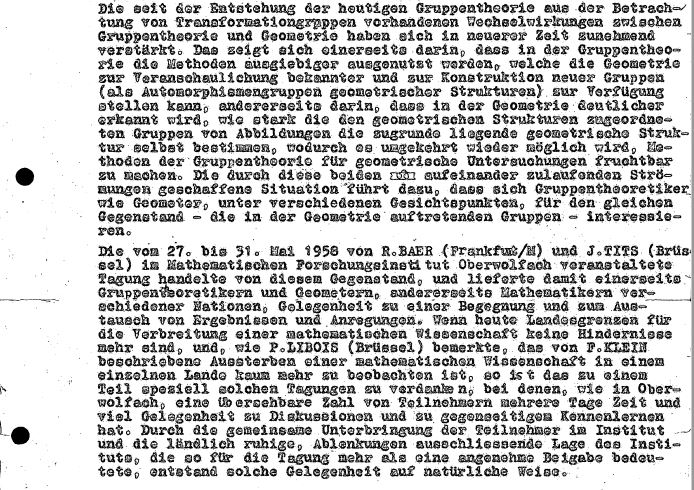}\caption{Auszug aus dem Tagungsbericht 1958.}

\end{figure}

\smallskip
Es ist deutlich erkennbar, dass sich der Bericht nicht nur auf fachliche Inhalte und auf die Arbeitsatmosph\"are bezieht, sondern dass er auch f\"ur eine wissenschaftspolitische Stellungnahme genutzt wird. Weiterhin sind die im Bericht formulierte \textit{\glqq \"ubersehbare Anzahl von Teilnehmern{\grqq}} und die Abgeschiedenheit wesentliche Bestandteile des \glqq Oberwolfacher Stils{\grqq}, der meines Wissens nach im Tagungsbericht 1959\footnote{Siehe \cite{ODA}, Workshop Nr. 5921, Seite 1.} zum ersten Mal explizit erw\"ahnt wird, sich in fr\"uheren Berichten aber schon andeutet.
\textit{\glqq Die diesj\"ahrige Tagung 1959 drohte durch das Ausma{\ss} der Beteiligung (40 aktive Teilnehmer und 25 Vortr\"age) den Rahmen des Oberwolfacher Kolloquiums zu sprengen. Das Interesse f\"ur diese Veranstaltung, die pers\"onliche Begegnungen mit den meisten prominenten Gruppentheoretikern Europas erm\"oglicht, ist inzwischen so gro{\ss} geworden, da{\ss} bei den weiteren Tagungen mit einer Teilung oder mit beschr\"ankter Teilnehmerzahl gerechnet werden mu{\ss}, um den \glqq Oberwolfacher Stil{\grqq} zu wahren, der durch eine sehr intensive Arbeitsatmosph\"are und einen sehr intensiven Gedankenaustausch au{\ss}erhalb der Vortr\"age gekennzeichnet ist. Den beiden wissenschaftlichen Tagungsleitern, Prof. R. Baer (Frankfurt a.M.) und Prof. H. Wielandt (T\"ubingen), ist es zu danken, da{\ss} durch eine sehr straffe Organisation und rigorose Redezeitbeschr\"ankungen der \glqq Oberwolfacher Stil{\grqq} bei dieser Tagung gewahrt bleiben konnte.{\grqq}}

Die Tagung 1959 hat laut Bericht den Anspruch, einen \glqq ziemlich repr\"asentativen Querschnitt{\grqq}\footnote{Ebenda.} durch aktuelle Forschungsthemen zu zeigen, inklusive der Pr\"asentation bisher unver\"offentlichter Ergebnisse von Helmut Wielandt und John G. Thompson.
Zur Tagung 1960 liegt leider kein Bericht vor, aber das Vortragsbuch (siehe \cite{ODA}) zeigt eine gro{\ss}e thematische Bandbreite und zahlreiche Vortragende aus dem Ausland (v.a. UK). Au{\ss}erdem ist zum ersten Mal explizit im Vortragsbuch erkennbar, dass Personen \"uber neue Resultate vorgetragen haben, die in der Zusammenarbeit mit anderen, teils ebenfalls anwesenden Personen erarbeitet worden waren\footnote{Gut zu erkennen auf den Seiten 219 und 221 des Vortragsbuchs Nr. 5.}. Dies k\"onnte ein Hinweis darauf sein, dass immer noch die wissenschaftspolitische Bedeutung des MFO unterstrichen werden sollte und dass hier ein starkes Argument gesehen wurde. Schlie{\ss}lich war erst kurz vorher, im Juli 1959, die Gesellschaft f\"ur Math. Forschung e.V. gegr\"undet worden, als Tr\"agerverein f\"ur das Forschungsinstitut.
Im Bericht zur Gruppentheorie-Tagung im Oktober 1961 wird betont, dass der Kreis der Eingeladenen etwas kleiner war, so dass es mehr Zeit f\"ur \glqq ausgedehnte Diskussionen{\grqq} gab, dass in den 20 Vortr\"agen trotzdem ein guter Querschnitt \textit{\glqq durch die heutigen Forschungsgebiete der Gruppentheorie{\grqq}} gegeben wurde und dass f\"unf amerikanische Mathematiker teilgenommen haben.
Diese sind namentlich im Bericht genannt: Michio Suzuki (Urbana), Donald G. Higman and Daniel Hughes (Ann Arbor), James S. Frame (East Lansing) and 
W. Holland (Tulane).
Generell f\"allt die Betonung der Teilnehmer*innen aus dem Ausland in den fr\"uhen Berichten auf. Erst im Lauf der 60er Jahre lesen sich die Berichte so, als sei das gemischte Publikum inzwischen Normalit\"at. Im Tagungsbericht zu \glqq Die Geometrie der Gruppen und die Gruppen der Geometrie
unter besonderer Ber\"ucksichtigung endlicher Strukturen{\grqq} aus dem Jahr 1964\footnote{Siehe Workshop Nr. 6421 in \cite{ODA}.} werden etwa die Teilnehmer*innen zwar nach L\"andern sortiert genannt, aber ohne besondere Betonung.
Da mir mehrere Zeitzeugen geschildert haben (und ich aus eigener Erfahrung best\"atigen kann), wie bereichernd die pers\"onlichen, informellen Gespr\"ache auf solchen Tagungen f\"ur sie waren und dass es eine der seltenen Gelegenheiten war, ber\"uhmte Kolleg*innen im Vortrag zu erleben, m\"ochte ich noch aus einem Bericht aus dem Jahr 1962 zitieren, der genau das thematisiert\footnote{Seite 1 des Berichts zur Tagung Gruppentheorie, 3.-10. August 1962, Workshop Nr. 6232 in \cite{ODA}}:
\textit{\glqq Die diesj\"ahrige Tagung \"uber Gruppentheorie in Oberwolfach unter Leitung von Professor R. BAER (Frankfurt a.M.) und Professor H. WIELANDT (T\"ubingen) gab durch die Anwesenheit einiger hervorragender Gruppentheoretiker aus dem Ausland Gelegenheit zu fruchtbarer Aussprache. Besonders wertvoll f\"ur die zahlreich anwesenden jungen deutschen Mathematiker war der pers\"onliche Kontakt mit diesen Wissenschaftlern.{\grqq}}
Ein kurzer Blick in die Anwesenheitsliste reicht, um diesen Einstieg in den Bericht nachvollziehen zu k\"onnen. Donald G. Higman, Noboru Ito, Michio Suzuki, John G. Thompson, Hans Zassenhaus, Jacques Tits, dann Reinhold Baer und Helmut Wielandt selbst - nat\"urlich hinterlie{\ss} ein solches Zusammentreffen bei den anwesenden Doktoranden einen tiefen Eindruck! Einige der jungen Wissenschaftler*innen von damals sollten sp\"ater selbst Ber\"uhmtheit erlangen. Besonders ein Anwesender, damals noch sehr jung, aber schon bekannt und als \textit{\glqq verdammt scharfsinniger Bursche{\grqq}}\footnote{Zitat von Wielandt, siehe Seite 327 in \cite{Solomon}.} gesch\"atzt, wird uns noch genauer besch\"aftigen.

Vorher schauen wir uns die eingeladenen Frauen n\"aher an.


\vspace{1cm}
\section{Frauen}

\textbf{Hel Braun} (Helene Braun, 1914-1986, Sch\"ulerin von Carl Ludwig Siegel) und \textbf{Paulette Libermann}\footnote{Photo von der Seite http://irma.math.unistra.fr/~maudin/Paulette-Photos.html} (1919-2007, Sch\"ulerin von Charles Ehresmann) nahmen 1954 teil, und Libermann hielt einen Vortrag \"uber Lie-Pseudogruppen und -Un\-ter\-grup\-pen.

\begin{figure}[h]

\includegraphics[scale=0.4]{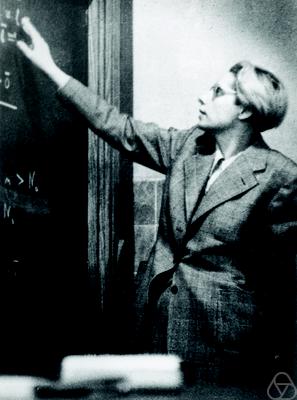}\caption{Hel Braun}

\end{figure}

Bei Hel Braun war es nicht der erste Besuch -- bereits im August und September 1951 ist sie im Vortragsbuch des MFO vermerkt, einmal mit
symmetrischen und hermitischen Matrizenpaaren und einmal mit hermiteschen Modulgruppen\footnote{Vortragsbuch 2, 4.4.1949-5.6.1952, auf der Seite \cite{ODA}.}. Diese Eintr\"age und auch die Teilnahme an der Tagung 1954 fallen in die 
Zeit, in der Braun Professorin an der Universit\"at Hamburg war und bevor sie den Lehrstuhl von Helmut Hasse \"ubernahm. 
Bei Paulette Liberman ist nur eine Tagungsteilnahme am MFO dokumentiert\footnote{Book of Abstracts 3, Zeitraum
1. Aug. 1952 -- 25. März 1955 auf der Seite \cite{ODA}.}, und aus den bekannten biographischen Angaben\footnote{Siehe zum Beispiel \cite{Libermann}.} 
geht nicht ganz klar hervor, ob sie noch in Stra{\ss}burg oder schon als Professorin in Rennes war.

\begin{figure}[h]

\includegraphics[scale=0.4]{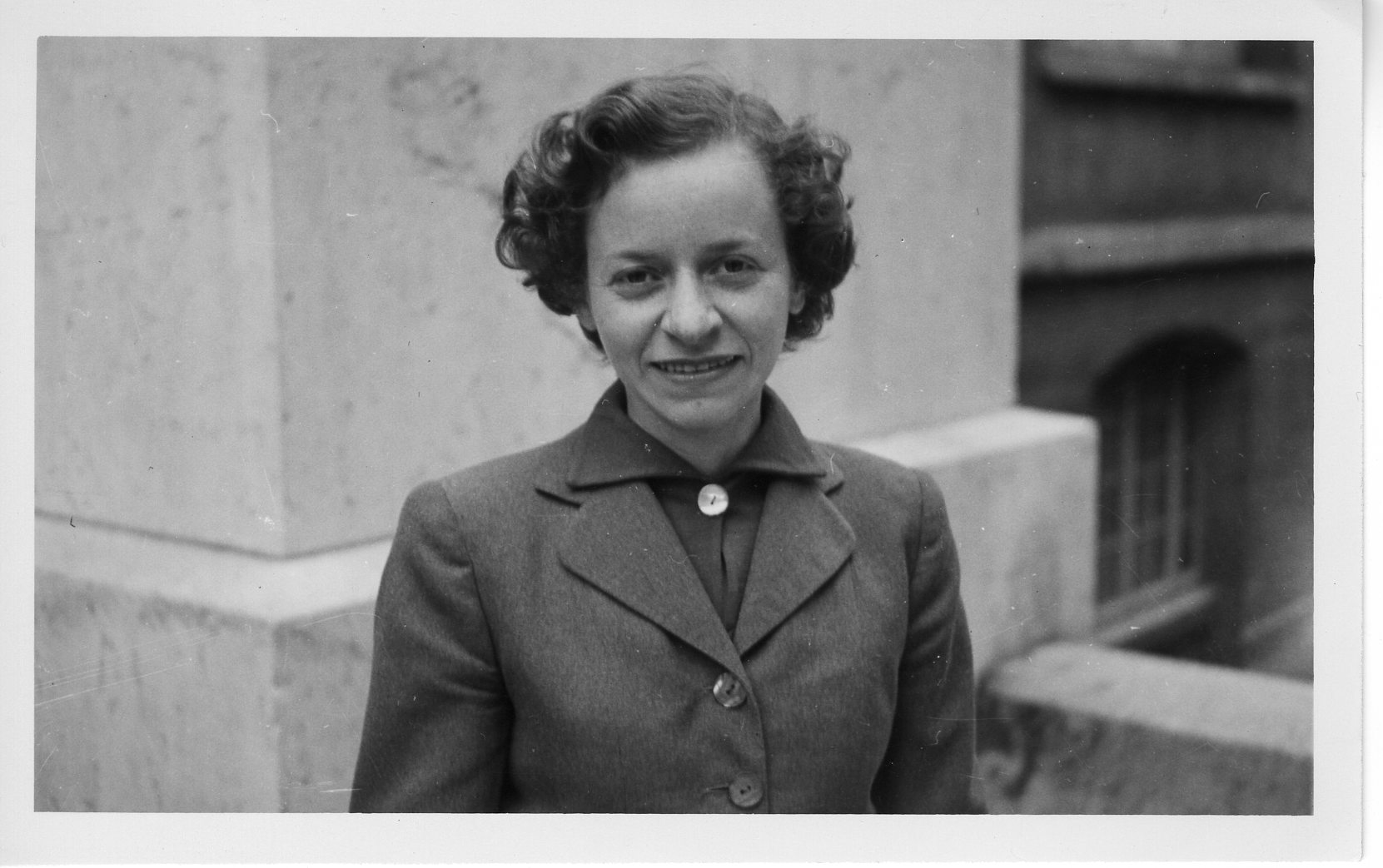}\caption{Paulette Liberman}

\end{figure}

\textbf{Olga Taussky-Todd} (1906-1995, geb. Tau{\ss}ky), eine Sch\"ulerin von Philipp Furtw\"angler, die aber auch stark von Emmy Noether beeinflusst war (siehe dazu \cite{Koreuber}), war 1955 dabei und trug vor \"uber Gruppenmatrizen. Zu der Zeit lebte sie noch in den USA und war nicht an einer Universit\"at t\"atig  (erst ab 1957 wieder). In der dazugeh\"origen Publikation (\cite{Taussky}) klassifiziert sie Gruppen, deren Gruppenmatrix gewisse Eigenschaften hat. Hier gibt es Zusammenh\"ange u.a. zu Dedekind-Gruppen, also zu Gruppen, in denen alle Untergruppen Normalteiler sind.

\begin{figure}[h]
\centering
\includegraphics{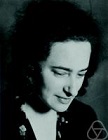}\caption{Olga Taussky-Todd}

\end{figure}

Taussky-Todd ist sehr produktiv, MathSciNet listet 155 Arbeiten mit insgesamt 503 Zitaten, und sie ist eine Pionierin bei algorithmischen Methoden, die etwa f\"ur die Entwicklung der Computeralgebra relevant sind. Aus den Unterlagen in\cite{ODA} geht hervor, dass sie sp\"ater noch mehrmals an Tagungen in Oberwolfach teilnahm, oft gemeinsam mit ihrem Mann John Todd und zu unterschiedlichen mathematischen Themen. Maths Genealogy (27.5.2021): \glqq Olga Taussky-Todd has 14 students and 161 descendants.{\grqq} Dabei kommen die meisten mathematischen Nachkommen \"uber Hanna Neumann zustande, die wir bereits erw\"ahnt haben.\\
 
\textbf{Maria Hasse} (1921-2014), eine Sch\"ulerin von Hans Schubert aus Rostock, ist 1958 dabei mit einem Vortrag zu Kategorien und Gruppoiden, und dann noch einmal 1959 mit S\"atzen \"uber Kategorien. Sie war bereits seit 1954 Professorin an der TU Dresden – die erste Professorin \"uberhaupt dort im Bereich Mathematik/Naturwissenschaften.

\begin{figure}[h]

\includegraphics[scale=0.4]{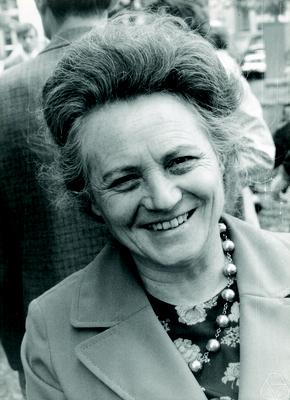}\caption{Maria Hasse}

\end{figure}

Im G\"astebuch findet sich ab 1958 mehrmals der Name \textbf{Griselda Pascual}. Ich vermute aufgrund des Arbeitsgebiets, dass es sich um die 1926 in Barcelona geborene und 2001 verstorbene Kollegin Pascual handelt, die unter der Betreuung von Enrique Linés Escardó promoviert wurde\footnote{Information in \cite{MGP}.}. Einige Informationen \"uber sie sind auf der Internetseite \cite{GP} zu finden: Bereits im Alter von 16 Jahren erlangte sie die Hochschulreife und legte alle n\"otigen Examina ab, um selbst an einer Schule unterrichten zu k\"onnen. Ihr Mathematikstudium schloss sie im Alter von 20 ab und arbeitete da bereits als Assistenzlehrkraft an der
Universitat de Barcelona, im Bereich der Geometrie.
1958 erhielt sie Stipendien vom Consejo Superior de Investigaciones Científicas (CSIC) und von der Humboldt-Stiftung und verbrachte u.a. Zeit in Freiburg, wo sie der Analysis zugeordnet war und an der Schnittstelle zur Algebra und Geometrie arbeitete.
Sie hatte immer wieder Positionen an Universit\"aten, wurde aber erst 1975 promoviert (mit der besten m\"oglichen Bewertung)
und wirkte sp\"ater bis zu ihrem Ruhestand als Professorin an der Universitat de Barcelona.
Zum Namen Anna Mahndler (G\"astebuch f\"ur 1961, schlecht zu lesen) konnte ich leider keine Angaben finden. Sp\"ater hat Baer noch einige Doktorandinnen, von denen aber nach meinen Informationen keine selbst wissenschaftlichen Nachwuchs hat. Sie sind also vermutlich nach der Promotion nicht (lange) weiter wissenschaftlich aktiv gewesen. Eine Frau f\"allt durch ihre Abwesenheit auf: \textbf{Ruth Moufang} (1905-1977), Sch\"ulerin von Max Dehn und die erste \"uberhaupt in Deutschland auf eine Mathematikprofessur berufene Frau, war zu der Zeit in Frankfurt aktiv und hat sp\"ater gemeinsam mit Baer zahlreiche Dissertationen betreut. Ihre Arbeit passt sehr gut zur Schnittstelle Gruppentheorie/Geometrie, weshalb es auf den ersten Blick erstaunlich ist, dass wir ihren Namen kein einziges Mal im Zusammenhang mit diesen fr\"uhen
Tagungen sehen.

\begin{figure}[h]

\includegraphics[scale=0.5]{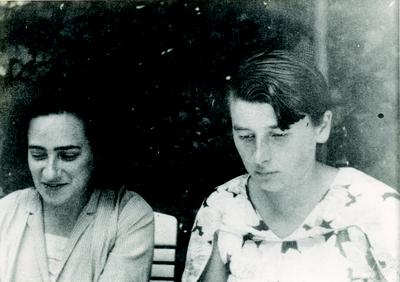}\caption{Olga Taussky-Todd und Ruth Moufang, 1932}

\end{figure}

Renate Tobies gibt in ihrem Einf\"uhrungskapitel in \cite{Tobies} detaillierte Einblicke in die vielf\"altigen Schwierigkeiten, die Frauen bei einer wissenschaftlichen Laufbahn begegnen\footnote{ab Seite 50, speziell zur Laufbahn als Hochschullehrerin ab Seite 53}.
So stellte etwa die Habilitation lange eine gro{\ss}e H\"urde da, und die ohnehin schon schwierige Situation f\"ur Wissenschaftlerinnen versch\"arfte sich ab 1933 weiter, nicht nur f\"ur J\"udinnen\footnote{Erster Absatz auf Seite 177 in \cite{Pieper}.}. Ruth Moufang wurde z.B. 1937 in einem Schreiben dar\"uber informiert, dass Frauen generell ungeeignet seien, Studenten, also m\"annliche Studierende, auszubilden. 
Max Dehn setzte sich zwar f\"ur Moufang ein, konnte das aber aus politischen Gr\"unden nicht effektiv tun (Amtsenthebung 1935).
In \cite{Pieper} sind detaillierte Hintergrundinformationen dazu zu finden -- zahlreiche Frauen verlie{\ss}en nach der Dissertation das wissenschaftliche Umfeld und wechselten in den Schuldienst oder in die Industrie. So auch Ruth Moufang, die jahrelang in Essen bei Krupp t\"atig war\footnote{Seite 187 in \cite{Pieper}} und parallel weiterhin wissenschaftlich arbeitete und publizierte. 
Lang war daher auch die Pause zwischen ihrer Habilitation (Frankfurt 1937) und der venia legendi (Frankfurt 1946), es folgte ein kommissarisch wahrgenommenes Extraordinariat ab 1947 und erst 1951 die Berufung nach Frankfurt als au{\ss}erordentliche Professorin.
Sie warb immer wieder Mittel f\"ur Gastprofessuren ein, und der erste Eingeladene war 1952/1953 Reinhold Baer, zu dem es eine jahrelange enge wissenschaftliche Verbindung gab mit zahlreichen gemeinsam betreuten Arbeiten (\cite{Pieper} ab Seite 191, auch \cite{MGP}).
Anders als Emmy Noether etablierte Ruth Moufang keine Schule. Ihre Erfahrungen, als Frau trotz gro{\ss}er wissenschaftlicher Leistungen \"uber viele Jahre nicht akzeptiert zu werden, hinterlie{\ss}en ein Gef\"uhl tiefer Kr\"ankung, das sie lange begleitet hat \footnote{Siehe \cite{Pieper}, Seite 194.}. Die wenigen Tagungsteilnahmen werden in \cite{Pieper} ebenfalls thematisiert\footnote{Seite 190 dort} und mit famili\"aren Umst\"anden in Zusammenhang gebracht. 
Interessant w\"are es, n\"aher zu untersuchen, wie sie ihr direktes mathematisches Umfeld in Frankfurt ab 1947 beeinflusst hat und inwiefern sich das nach der Ankunft von Reinhold Baer ver\"andert hat.

\section{Thematische Schlaglichter}

Hier geht es um einige Personen und deren Tagungsbeitr\"age, die noch nicht besonders beleuchtet wurden. Die Auswahl erfolgte komplett willk\"urlich: In zwei F\"allen wurde \"uber Resultate vorgetragen, die in einer Geschwindigkeit bekannt und verbreitet waren, die aus heutiger Sicht erstaunlich ist, und die bis heute ganz klar mit den Urhebern verkn\"upft sind. In den beiden anderen F\"allen handelt es sich um Pers\"onlichkeiten, die die weitere Entwicklung der Gruppentheorie entscheidend mit gepr\"agt haben.\\

\textbf{Otto Gr\"un} (1888-1974) ist u.a. deshalb ein interessanter Fall, weil er nicht im klassischen Sinne, also mit einer Assistentenstelle oder Professur, im akademischen Umfeld etabliert war. Laut Wikipedia  war er Autodidakt, korrespondierte u.a. mit Helmut Hasse und ver\"offentlichte zahlreiche mathematische Arbeiten mit Beitr\"agen zur Theorie der endlichen Gruppen. Hier ein Beispiel (siehe \cite{Gruen1960}), das auf umfangreichen Vorarbeiten des Autors und relevanten Resultaten von Helmut Wielandt beruht: \glqq F\"ur jede Primzahl p ist die maximale p-Faktorgruppe von G eindeutig bestimmt. Die maximale nilpotente Faktorgruppe von G ist eindeutig bestimmt und das direkte Produkt der maximalen p-Faktorgruppen (f\"ur die relevanten Primzahlen p).{\grqq} Die nach ihm benannten S\"atze sind heute ganz selbstverst\"andliches Lehrbuchmaterial in der Gruppentheorie. Er besuchte mehrere Tagungen in Oberwolfach w\"ahrend seiner Zeit als Lehrbeauftragter an der Universit\"at W\"urzburg und hielt regelm\"a{\ss}ig Vortr\"age, war also durchaus in die mathematische Community eingebunden. 1955 trug er vor \"uber freie Gruppen (Morphismen von Faktorgruppen etc.), 1959 \"uber p-Faktorgruppen (siehe \cite{Gruen1960}) und 1960 \"uber m-aufl\"osbare und m-nilpotente Gruppen. Den Vortrag 1961 mit dem Titel \glqq Im Gruppenring von G konjugierte Untergruppen von Gasch\"utz{\grqq} konnte ich trotz Literaturrecherche bisher nicht einordnen.\\

\textbf{Roger William Carter} (geb. 1934) kam vermutlich aus Cambridge zur Tagung 1959, da er 1960 in Cambridge unter der Betreuung von Derek Roy Taunt promoviert wurde. Er trug bei dieser Tagung zu Sylowsystemen in aufl\"osbaren Gruppen vor. 1960 k\"onnte er aus T\"ubingen angereist sein, wo er in etwa diesem Zeitraum eine Humboldt Research Scholarship innehatte\footnote{Siehe Fu{\ss}note auf der ersten Seite von \cite{Carter1961}, wurde im September 1960 eingereicht.}. Dauer, Zeitpunkt und genaue Umst\"ande des Gastaufenthaltes konnte ich bisher nicht herausfinden. Carter arbeitete zwischenzeitlich viel zu aufl\"osbaren Gruppen und ist bis heute bekannt f\"ur die jetzt nach ihm benannten Gruppen, \"uber die er in seinem Vortrag 1960 sprach. Das Hauptresultat der dazugeh\"origen Publikation ist: \textit{Jede endliche aufl\"osbare Gruppe hat nilpotente Untergruppen, die ihr eigener Normalisator in der Gruppe sind, und all diese sind konjugiert.} Carter schrieb dazu in \cite{Carter1961}: \textit{\glqq I would like to thank Professor WIELANDT and Professor P. HALL for their help in the formulation of the main theorem in its present general form and in the construction of the proof.{\grqq}} Aufl\"osbare Gruppen waren \"uber Jahrzehnte ein wichtiger Forschungsgegenstand, nicht nur wegen der Beziehungen zur Galoistheorie. Carters Arbeit wirft viele Folgefragen f\"ur nicht-aufl\"osbare Gruppen auf, und das f\"uhrt u.a. zu Charakterisierungen gewisser einfacher nicht-abelscher Gruppen. Es gibt die Vermutung, dass die heute Carter-Untergruppen genannten Gruppen dann, wenn sie existieren, immer konjugiert sind. Daran wird immer noch geforscht, was ein Grund f\"ur die zahlreichen Zitate einiger seiner Arbeiten ist.\footnote{ MathSciNet listet 53 Publikationen und 2943 Zitate.} Einflussreich waren auch seine Beitr\"age zur 2-Sylowstruktur der klassischen Gruppen, zusammen mit Paul Fong, und es gibt einige Standardb\"ucher von ihm zu endlichen Gruppen vom Lie-Typ. 

\textbf{Joachim Neub\"user} (1932--2021) war ein Sch\"uler von Wolfgang Gasch\"utz und arbeitete auch nach der Promotion noch in Kiel, bevor er 1969 an der RWTH Aachen den in Fachkreisen bekannten, inzwischen umbenannten Lehrstuhl D \"ubernahm. Auf einer Zuse 22 in Kiel am Rechenzentrum (eingerichtet von Karl Heinrich Weise) wurden um 1960 herum Anf\"ange der computergest\"utzten Gruppentheorie umgesetzt. Neub\"user beschreibt in \cite{Neubuser} die theoretischen Grundlagen f\"ur die technische Umsetzung, also daf\"ur, dass Gruppenelemente \"uberhaupt im Computer repr\"asentiert werden k\"onnen und dass mit ihnen gearbeitet werden kann. Was wir inzwischen ganz selbstverst\"andlich als Computeralgebra kennen, hei{\ss}t bei ihm \glqq Untersuchungen auf einer programmgesteuerten elektronischen Dualmaschine{\grqq}.
Am Ende seines Artikels \cite{Neubuser} beschreibt er die technischen Details der Z 22, gibt Rechenzeiten an, und sowohl er selbst als auch sein Kollege Volkmar Felsch k\"onnen bzw. konnten lebhaft von diesen Anf\"angen der Computeralgebra berichten\footnote{Telefonate mit beiden im M\"arz 2021.} – inklusive der Arbeit mit Lochstreifen zum Einlesen der Programme. Neub\"user trug 1959 und 1960 in Oberwolfach \"uber die ersten Ergebnisse vor und motivierte seine Pionierarbeit: Der Computer liefert Beispielmaterial f\"ur Vermutungen etc. und hilft dort weiter, wo noch wenige systematische Methoden zur Strukturuntersuchung vorhanden sind. Die abstrakte Gruppentheorie war immerhin damals noch in der Entstehungsphase! Beispiele f\"ur den Einsatz von Computern: Eine Gruppe ist gegeben durch Erzeugende und Relationen, dann erfolgt z.B. eine Restklassenabz\"ahlung f\"ur eine Permutationsdarstellung. Rechenverfahren sollen m\"oglichst auf verschiedene Klassen von Gruppen anwendbar sein, und innerhab eines Algorithmus werden dann gegebenenfalls Verzweigungen entlang von Spezialf\"allen eingebaut. Neub\"user berechnet in \cite{Neubuser}  bis zu einer gewissen Gruppengr\"o{\ss}e den Untergruppenverband (den sogenannten \glqq Situationsplan{\grqq}), und darin sind gewisse charakteristische Untergruppen ablesbar. 

\begin{figure}[h]

\includegraphics[scale=0.6]{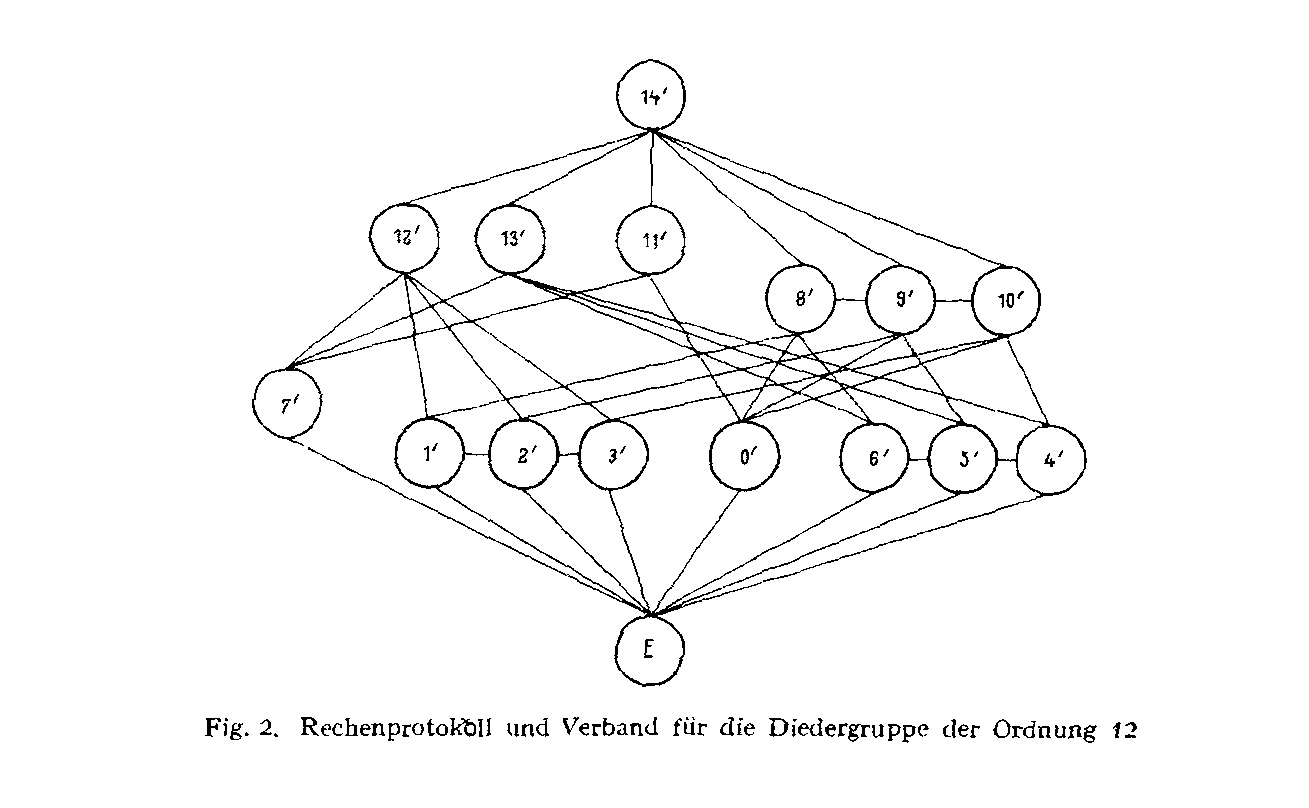}\caption{Abbildung des Untergruppenverbandes von $D_{12}$ aus \cite{Neubuser}}

\end{figure}

Das Rechenverfahren wird erkl\"art, wobei der Artikel sich auf Permutationsgruppen konzentriert. Wesentliche Unterprogramme, die beschrieben werden: Arithmetisches, z.B. Rechnen mit ganzen Zahlen bis hin zur Primfaktorzerlegung, Lesen, Drucken, Umschreiben, Vergleichen. Intern wird mit Abbildungen gerechnet, f\"ur die optische Darstellung von Permutationen wird die Zyklenschreibweise verwendet, und das Programm kann bereits verschiedene Rechnungen mit Permutationen durchf\"uhren. Nur wenige Grundrechenoperationen sind fest in der Maschine verdrahtet – der Rest wird als Befehl mit Lochstreifen eingegeben. Neub\"users Artikel  zeigt Beispiele f\"ur die vom Programm ausgegebenen Informationen und gibt Rechenzeiten an. Seine fr\"uhen Arbeiten werden laut \cite{MSN} bis jenseits von 2000 zitiert, u.a. f\"ur Markentafeln, die ihrerseits auch wieder zahlreiche Anwendungen haben. Joachim Neub\"user wurde viel sp\"ater selbst Tagungsleiter bei den Oberwolfach-Workshops zur Computeralgebra, und vorher gab es bereits Tagungen zu \glqq Grundlagen des numerischen Rechnens und Computeralgebra{\grqq}. 
Mehrere seiner Sch\"uler*innen sind noch an Universit\"aten aktiv, und zahlreiche mathematische Nachkommen gehen auf Wilhelm Plesken zur\"uck\footnote{Maths Genealogy am 27.5.2021: \glqq Joachim Neub\"user has 9 students and 59 descendants.{\grqq}} MathSciNet listet 28 Publikationen und 212 Zitate, wovon mit Abstand die meisten auf das Buch \cite{BBNWZ} verweisen.\\

Der damals noch sehr junge \textbf{John Griggs Thompson} (geb. 1932), ein Sch\"uler von Saunders Mac Lane, war f\"ur inhaltliche H\"ohepunkte auf mehreren Tagungen verantwortlich. Er ist der \glqq verdammt scharfsinnige junge Bursche{\grqq}, der zu Beginn erw\"ahnt wurde. W\"ahrend er an einem Beweis der Vermutung arbeitete, dass Gruppen ungerader Ordnung aufl\"osbar sind\footnote{Bewiesen mit Walter Feit in \cite{FT1963}.}, trug er mehrmals in Oberwolfach vor. Gerhard Pazderski war dabei, als Thompson bei der Tagung 1960 \"uber den wichtigen Spezialfall berichtete, dass in einem minimalen Gegenbeispiel zum \glqq Odd order theorem{\grqq} die Sylowuntergruppen abelsch sind. Ich zitiere hier aus meiner schriftlichen Kommunikation mit ihm\footnote{Kommunikation per E-Mail am 5. Mai 2021.}: \textit{\glqq Eine herausragende Rolle spielte J. G. Thompson, den ich bei meiner letzten Vorwendetagung 1960 erlebte. Er hatte kurz zuvor seine aufsehenerregenden Ergebnisse \"uber fixpunktfreie Automorphismen im Zusammenhang mit der Nilpotenz des Frobeniuskerns gefunden und dar\"uber vorgetragen. In Arbeit hatte er zusammen mit Feit den Odd Order-Satz und inzwischen schon Spezialf\"alle erledigt, wenn ich mich recht erinnere den mit der Zusatzvoraussetzung abelscher Sylowgruppen. Er wurde gebeten im Verlauf der Tagung etwas \"uber diese Untersuchungen zu erz\"ahlen, was er auch bereitwillig tat. Was hier ablief, war f\"ur mich eine Sternstunde der Gruppentheorie.{\grqq}} 
W\"ahrend es bei diesem Vortrag einen Eintrag im Vortragsbuch gibt und ein Abstract, das zumindest f\"ur Eingeweihte den spektakul\"aren Inhalt dieses Tagungsbeitrags deutlich macht, muss man bei der n\"achsten Sensation in den dazugeh\"origen Tagungsbericht schauen. John G. Thompson erscheint n\"amlich 1962 weder im G\"aste- noch im Vortragsbuch, aber aus dem Bericht geht hervor, dass die \glqq Odd order-Arbeit{\grqq} jeden Nachmittag diskutiert wurde. Genauer hei{\ss}t es dort\footnote{Seite 1 des Berichts zur Tagung Gruppentheorie, 3.-10. August 1962, Workshop Nr. 6232 in \cite{ODA}.} dazu: \textit{\glqq Im Mittelpunkt der Tagung standen die an jedem Nachmittag gehaltenen
Vortr\"age von J. THOMPSON \"uber den von ihm gemeinsam mit W. FEIT bewiesenen Satz von der Aufl\"osbarkeit der Gruppen ungerader Ordnung. Der sehr schwierige Beweis dieses seit langem vermuteten Satzes ist noch nicht ver\"offentlicht und wird auch in n\"achster Zeit, schon wegen der au{\ss}erordentlichen L\"ange von 400 Seiten Schreibmaschinenschrift, nicht ver\"offentlicht werden. Es war daher f\"ur die Teilnehmer sehr n\"utzlich, die Methoden kennenlernen zu k\"onnen, denn nur eine Tagung im \glqq Oberwolfacher Stil{\grqq} erlaubt eine ausf\"uhrliche Diskussion einer solchen umfangreichen Arbeit.{\grqq}}

\begin{figure}[h]

\includegraphics{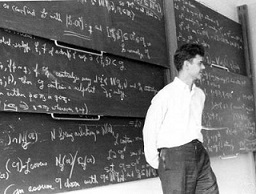}\caption{John G. Thompson}

\end{figure}

Da ist er wieder, der \textbf{Oberwolfacher Stil}.
MathSciNet listet f\"ur John Griggs Thompson 99 Publikationen und 1307 Zitate, und unter den die am meisten zitierten Arbeiten sind gleich zwei\footnote{n\"amlich \cite{Thompson} und \cite{FT1963}}, \"uber die er im hier betrachteten Zeitraum auf Tagungen in Oberwolfach vorgetragen hat. Das Pazderski-Zitat oben weist n\"amlich auch auf einen Vortrag 1960 hin, der das Thema von Thompsons Dissertation behandelt. Ron Solomon schreibt dazu in \cite{Solomon}\footnote{auf Seite 324, Mitte} in einem Abschnitt, der passenderweise mit \glqq Enter John Thompson{\grqq} beginnt: \textit{\glqq At the suggestion of Marshall Hall, Thompson attacked in his dissertation the long-standing conjecture that the Frobenius kernel is always nilpotent.{\grqq}}
Kurz danach\footnote{unten auf der gleichen Seite} hei{\ss}t es weiter: \textit{\glqq Thompson's thesis had immediate implications for the study of Zassenhaus groups, clarifying as it did the structure of Frobenius kernels. Even more important, it was the beginning of Thompson's profound analysis of the structure of solvable subgroups of finite simple groups. In the summer of 1958 while working on his thesis, Thompson visited Wielandt in T\"ubingen, and Huppert reports Wielandt's comment: Das ist ein verdammt scharfsinniger Bursche. Man kann etwas lernen von ihm. [That's one damn sharp guy. You can learn something from him.]{\grqq}}
Wir erleben hier in Oberwolfach also mit, wie Thompson damit beginnt, durch seine Ideen die weitere Entwicklung der endlichen Gruppentheorie ma{\ss}geblich zu beeinflussen. Dass hier ein Forschungsgebiet zu neuem Leben erwacht, zeigt sich auch an anderer Stelle: \textit{\glqq By 1959 when Marshall Hall published his text \glqq The Theory of Groups{\grqq}, he could write in dramatic contrast to Brauer's remarks in 1954: \glqq Current research in Group Theory, as witnessed by the publications covered in Mathematical Reviews, is vigorous and extensive.{\grqq}} Solomon bezieht sich hier in \cite{Solomon} auf Richard Brauers Bemerkungen 1954 beim ICM, dass es eine gewisse Stagnation in der Entwicklung der Theorie der endlichen Gruppen gebe. Die Entwicklungen hin zum Klassifikationsprojekt zeigen, dass sp\"atestens mit diesen wichtigen Arbeiten von Thompson um 1960 herum die Stagnation ein Ende hatte.

\section{Netzwerke und beginnende Traditionen}

Wir beleuchten nun noch einmal die Schl\"usselfiguren und deren Standorte zur Zeit der hier behandelten Tagungen.
Was die Planung der Gruppentheorie-Tagungen betrifft, f\"allt besonders die Kommunikation zwischen Wolfgang Gasch\"utz, Wilhelm S\"uss und Helmut Wielandt 
auf, und die stabilen Verbindungen ins Ausland, dabei mehr in Richtung Westen, sehen wir geb\"undelt in Kontaktpersonen, \"uber die die Korrespondenz l\"auft. S\"uss hat bewusst Kontakte zu emigrierten Kollegen gepflegt, was z.B. von Volker Remmert in \cite{Remmert}\footnote{auf Seite 361 dort} folgenderma{\ss}en eingeordnet wird: \textit{\glqq Suss was also interested in bringing Reinhold Baer, Friedrich Wilhelm Levi and Bernhard Neumann to the institute. Their reactions to his invitations were positive. Baer and Neumann would certainly have come to Oberwolfach as early as the late 1940s but for lack of travel funds. Levi came to Oberwolfach in 1950 and Baer in 1952 (not having managed to incorporate the visit into his 1950 travel plans). In 1951 Bernhard Neumann (1909-2002), who had emigrated to Great Britain in 1933, came with his wife Hanna (1914-1971), a group theorist as well. For Baer and Levi, these early visits to Oberwolfach were important steps on their way back to Germany.{\grqq}}

Eine Verbindungsfigur ist dabei Alfred Loewy, der 1919 bis 1934 Professor in Freiburg war
 und dem Wilhelm S\"uss auf die Professur in Freiburg nachfolgte.
Eine direkte Verbindung nach Kiel wird sichtbar in Karl-Heinrich Weise, dem Doktorvater von Wolfgang Gasch\"utz, der 
gemeinsam mit Reinhold Baer
zu den Gr\"undungsmitgliedern der Gesellschaft f\"ur Mathematische Forschung e.V. geh\"ort, dem 1959 gegr\"undeten Tr\"agerverein des MFO
\footnote{Siehe https://www.mfo.de/about-the-institute/structure/society/founder-members.}.
Auch Hellmuth Kneser finden wir in der Liste der Gr\"undungsmitglieder -- dieser hatte 
sich u.a. daf\"ur eingesetzt, dass Helmut Wielandt eine Assistentenstelle in Tübingen bekam, 
und er leitete das MFO kurzzeitig nach Wilhelm S\"uss. 
Ich vermute, dass die gr\"o{\ss}ere finanzielle Stabilit\"at f\"ur den Betrieb des MFO und der noch
st\"arkere Einfluss von Reinhold Baer nach seinem Ruf nach Frankfurt mit ausschlaggebend daf\"ur waren, wie
sich die Tagungstradition am MFO in den 1960er Jahren (zumindest in der Gruppentheorie) 
entwickelte, und es sprengt den Rahmen dieses Beitrags, das 
genauer zu beschreiben. Wir haben ja bereits gesehen, dass um 1960 herum auch inhaltlich pl\"otzlich 
sehr viel Spannendes in der Theorie der endlichen Gruppen passierte, so dass sich 
eine ausf\"uhrliche Befassung mit der weiteren Tagungstradition, der Publikationskultur und der Kommunikation der
nun in rasantem Tempo folgenden neuen Resultate lohnt.

Hier beschr\"anke ich mich auf weitere Bemerkungen zur Organisation dieser fr\"uhen Tagungen und zu den Traditionen, deren Entstehen wir beobachten k\"onnen. Wir haben Beispiele gesehen f\"ur die inhaltlichen Absprachen und die Einladungslisten -- ob bei den Personen auf den Einladungslisten das aktuelle Arbeitsgebiet ganz genau zur Ausrichtung der Tagung passt, schien zweitrangig zu sein, jedenfalls zeigt sich in der Korrespondenz nicht, dass Einladungen standardm\"a{\ss}ig davon abh\"angig waren, ob jemand gerade einen besonders spannenden Satz bewiesen hatte. 
Eher deuten die Votragsthemen und die inhaltlich passenden Publikationen darauf hin, dass es oft umgekehrt war und dass
der Vortrag in Oberwolfach als eine der ersten Kommunikationsm\"oglichkeiten f\"ur neue Resultate genutzt wurde, auch f\"ur 
Verfeinerungen vor der Publikation. Die Tagungen waren inhaltlich breit geplant und es ist nur selten erkennbar (z.B. bei der Korrespondenz mit Wolfgang Gasch\"utz im Sommer 1955), dass jemand \"uberhaupt konkrete Vortragsthemen vorher vorgeschlagen hat. Auch konnte ich keine Hinweise darauf finden, dass Personen aus konkreten Gr\"unden nicht wieder eingeladen wurden. Ein genauerer Blick in die G\"asteb\"ucher zeigt, dass die Organisatoren h\"aufig ihre Familie mit nach Oberwolfach brachten (siehe zum Beispiel in Abbildung 1, wo wir \glqq Irmgard Wielandt mit Eltern{\grqq} lesen). Wir sehen auch, dass sich fr\"uh etabliert hat, was wir auch von heutigen Fachtagungen kennen: Wer hinf\"ahrt, nimmt oft mehrere Doktorand*innen mit, und dem Nachwuchs wird m\"oglichst auch die Gelegenheit zum Vortrag gegeben. Dazu m\"ochte ich noch einmal aus der Korrespondenz mit Gerhard Pazderski\footnote{Kommunikation per E-Mail am 5. Mai 2021.} zitieren:
\textit{\glqq Ich war als junger Assistent mit meiner Dissertation besch\"aftigt, die ich mich bem\"uhte ohne meinen 1954 verstorbenen Doktorvater Heinrich Brandt auf dem Gebiet der Gruppentheorie anzufertigen. Obwohl kein Gruppentheoretiker war der hallesche Geometer Ott-Heinrich Keller, der auch \"uber Gruppen publiziert hat, bereit nach Brandts Tod meine Betreuung zu \"ubernehmen. Er nutzte seine Bekanntschaft mit Helmut Wielandt, dem damals f\"uhrenden Gruppentheoretiker in Deutschland, um mir den Besuch der j\"ahrlich in Oberwolfach durchgef\"uhrten Fachtagungen zu erm\"oglichen. Bis zum Bau der DDR-Grenzbefestigungen 1961 konnte ich drei Mal an den Zusammenk\"unften teilnehmen, zuletzt 1960. Diese Tagungsbesuche wirkten auf mich ungeheuer stimulierend.{\grqq}}

Zu Beginn der 1960er Jahre, also gegen Ende des hier betrachteten Zeitraums, wird Reinhold Baer noch mehr zur pr\"agenden Figur, mit zahlreichen kurzen Arbeitstagungen des Frankfurter Seminars und regelm\"a{\ss}igen \glqq Kindertagungen{\grqq}. Der Name der Tagung \glqq B\"aren mit Kindern und Kegeln{\grqq} (mit Bezug auf Baers Doktoranden Otto Kegel) ist tats\"achlich genau so im G\"astebuch zu finden. Im Januar 1963 gibt es dann den \glqq Baerschen Kinderausflug{\grqq}, und auch Zeitzeugen erinnern sich daran, dass stets in der ersten Januarwoche die Baersche Kindertagung dran war. Mehrere gr\"o{\ss}ere Tagungen werden zusammen mit Helmut Wielandt organisiert, oft mehrere innerhalb eines Jahres. Baers Engagement und die zahlreich stattfindenden Tagungen wirken sehr belebend auf die Entwicklung der Gruppentheorie, laut Zeitzeugen gerade zur richtigen Zeit: Mit John G. Thompson und seinen Arbeiten erfuhr die Entwicklung der abstrakten Theorie der endlichen Gruppen einen enormen Schub, und es entstand das, was wir heute \glqq lokale Theorie{\grqq} nennen. Die endliche Gruppentheorie rechnet in \glqq vor und nach Thompson{\grqq}, wie auch Solomon in \cite{Solomon}\footnote{ab Seite 324} verdeutlicht. Sein oben erw\"ahnter, gemeinsam mit Walter Feit bewiesener Satz, dass Gruppen ungerader Ordnung aufl\"osbar sind, sorgt f\"ur Schwung im Klassifikationsprojekt, denn nun ist endlich klar, dass jede einfache nicht-abelsche Gruppe ein Element der Ordnung 2 besitzt. Laut Brauer-Fowler (\cite{BF}) ergibt sich so eine Strategie, um entlang des Zentralisators eines solchen Elements die endlich vielen M\"oglichkeiten f\"ur die Struktur der ganzen Gruppe zu klassifizieren. Der internationale Besuch in Oberwolfach ist f\"ur viele Forschende die einzige Gelegenheit, mit Kolleg*innen aus dem Ausland und deren aktueller, teils noch unver\"offentlichter Arbeit in Kontakt zu kommen und umso besser mitverfolgen zu k\"onnen, wie neuer Schwung in das Gebiet kommt und ein L\"andergrenzen sowie Generationen \"uberschreitendes Momentum hin zum Klassifikationsprojekt entsteht.
Die Gruppentheorie verzweigt sich gleichzeitig immer weiter. Neben Tagungen im Format \glqq Gruppen und Geometrien{\grqq}, die es weiterhin gibt (zuletzt in Oberwolfach 2008, danach mehrmals an der BIRS Banff, Kanada), etablieren sich auch Tagungen zu Darstellungstheorie, Permutationsgruppen, zu algebraischen Gruppen, topologischen Gruppen und zu algorithmischer Gruppentheorie bzw. Computeralgebra.
Schlie{\ss}en m\"ochte ich mit einem Zitat aus dem Bericht zur Arbeitstagung des ehemaligen Baerschen Seminars 1986, geleitet von Helmut Salzmann\footnote{Workshop Nr. 8651, der Tagungsbericht ist in \cite{ODA} zu finden.}:
\textit{\glqq Das gro{\ss}e Interesse an dieser Tagung wird durch die Zahl von 41 Teilnehmern aus dem In- und Ausland dokumentiert, von denen 38 Vortr\"age hielten. Trotz der breit gestreuten Themen aus Algebra und Geometrie erwies sich die gemeinsame Wurzel und die pr\"agende Kraft von Baers mathematischer Denkweise als stark genug, ein inhaltliches Auseinanderfallen der Tagung zu verhindern. Im Gegenteil ergaben sich, vor allem auch in den Diskussionen, unerwartete Beziehungen zwischen entfernteren Gegenst\"anden. Dadurch war die Tagung anregender als viele andere, thematisch st\"arker eingeengte. Trotz des umfangreichen Vortragsprogramms bot die Atmosph\"are in Oberwolfach zahlreiche Gelegenheiten zum Erfahrungsaustausch und zu Problemsitzungen im kleinen Kreis. Herr B. Fischer (Bielefeld) konnte w\"ahrend der Tagung seinen 50. Geburtstag feiern.{\grqq}}

\vspace{0.5cm}
\textbf{Danksagung.}
Zahlreiche Personen haben in pers\"onlichen Gespr\"achen, am Telefon oder per E-Mail Erinnerungen an ihre ersten Oberwolfach-Tagungen mit mir geteilt, und daf\"ur bin ich sehr dankbar. Konkret hier eingeflossen sind Kommentare von Gerhard Pazderski, von Volkmar Felsch und von Joachim Neub\"user, der in der Zwischenzeit leider verstorben ist. Meinem Kollegen Ferran Dachs-Cadefau danke ich f\"ur wertvolle Hinweise zum Leben und Wirken von Griselda Pascual.


\end{document}